\documentclass[12pt,a4paper]{article}
\usepackage[latin2]{inputenc}
\usepackage{amsmath}
\usepackage{amsfonts}
\usepackage{amssymb}
\usepackage{graphicx}
\usepackage[left=2.00cm, right=2.00cm, top=2.50cm, bottom=2.50cm]{geometry}
\usepackage[T1]{fontenc}

\def\Q{{\mathbb Q}}
\def\Z{{\mathbb Z}}

\newtheorem{lemma}{Lemma}
\newtheorem{theorem}[lemma]{Theorem}
\newtheorem{corollary}[lemma]{Corollary}

\title{
Simplest quartic and simplest sextic Thue equations \\
over imaginary quadratic fields 
}
\author{
Istv\'{a}n Ga\'{a}l\thanks{
        Research supported in part by K115479 from the
        Hungarian National Foundation for Scientific Research				
				and by the EFOP-3.6.1-16-2016-00022 project. 
				The project is co-financed by the European Union and the European Social Fund.	
                         },\; \\
{\small University of Debrecen, Mathematical Institute} \\
{\small H--4002 Debrecen Pf.400., Hungary,} 
{\small e--mail: gaal.istvan@unideb.hu},
 \\ \\
Borka Jadrijevi\'c\thanks{
Research supported in part
 by the Croatian Science Foundation under the
project no. 6422.   }\\
{\small University of Split,}
{\small Faculty of Science, }\\
{\small Ru\dj era Bo\v{s}kovi\'{c}a 33, 21000 Split, Croatia,}
{\small e--mail: borka@pmfst.hr}
 \\ \\ 
L\'aszl\'o Remete\thanks{
        Research supported by the \'UNKP-17-3 new national excellence program of the ministry of human capacities.}\; \\
{\small University of Debrecen, Mathematical Institute} \\
{\small H--4002 Debrecen Pf.400., Hungary,} 
{\small e--mail: remete.laszlo@science.unideb.hu}
}

\begin{document}
\baselineskip=17pt

\maketitle
\thispagestyle{empty}

\renewcommand{\thefootnote}{\arabic{footnote}}
\setcounter{footnote}{0}

\noindent
Mathematics Subject Classification: Primary 11D59; Secondary 11D57\\
Key words and phrases: relative Thue equations, simplest quartic fields, simplest sextic fields

\begin{abstract}
The families of simplest cubic, simplest quartic and simplest
sextic fields and the related Thue equations are well known, see \cite{lpv}, \cite{lpv1}.
The family of simplest cubic Thue equations was already studied in the relative case,
over imaginary quadratic fields. In the present paper we give a similar extension
of simplest quartic and simplest sextic Thue equations over imaginary quadratic fields.
We explicitly give the solutions of these infinite parametric families of 
Thue equations over arbitrary imaginary quadratic fields.
\end{abstract}

\newpage

\section{Introduction}

Let $t$ be an integer parameter.
The infinite parametric families of number fields generated by the roots of the polynomials
\[
\begin{array}{lll}
&f_t^{(3)}(x)=x^3-(t-1)x^2-(t+2)x-1, &\;\;\;\; (t\in\Z),\\
&f_t^{(4)}(x)=x^4-tx^3-6x^2+tx+1,&\;\;\;\;(t\in\Z\setminus \{-3,0,3\}),\\
&f_t^{(6)}(x)=x^6-2tx^5-(5t+15)x^4-20x^3+5tx^2+(2t+6)x+1,&\;\;\;\;(t\in\Z\setminus \{-8,-3,0,5\}),
\end{array}
\]
are called {\bf simplest cubic, simplest quartic and simplest sextic fields}, respectively.
They are extensively studied in algebraic number theory, starting with D.Shanks \cite{shanks},
in the cubic case. 
It was shown by G.Lettl, A.Peth\H o and P.Voutier \cite{lpv} that these 
are all parametric families of number fields which are  
totally real cyclic with Galois group
generated by a mapping of type $x\mapsto \frac{ax+b}{cx+d}$ with $a,b,c,d\in\Z$.

Let $F(x,y)\in\Z[x,y]$ be an irreducible binary form of degree $\geq 3$
and let $0\neq k\in\Z$.
There is an extensive literature of {\bf Thue equations} of type 
\[
F(x,y)=k\;\;{\rm in}\;\;x,y\in\Z.
\]
In 1909 A.Thue \cite{thue} proved that these equations admit only finitely many solutions.
In 1967 A.Baker \cite{baker} gave effective upper bounds for the solutions.
Later on authors constructed numerical methods to reduce the bounds and 
to explicitly calculate the solutions, see \cite{book} for a summary.

In 1990 E.Thomas \cite{thomas} 
considered an {\bf infinite parametric family of Thue equations},
corresponding to the simplest cubic fields.
For $t\in\Z$, let
\[
F_t^{(3)}(x,y)=x^3-(t-1)x^2y-(t+2)xy^2-y^3
\]
and consider 
\[
F_t^{(3)}(x,y)=\pm 1\;\;{\rm in}\;\;x,y\in\Z.
\]
E.Thomas described the solutions for large enough parameters $t$,
later on the solutions were found for all parameters by M.Mignotte \cite{mig}.
This was the
first infinite parametric family of Thue equations that was completely
solved. Instead of single equations the solutions were given for 
infinitely many equations, for all values of the parameter $t$. 
These equations are also called 
the infinite parametric family of the 
{\bf simplest cubic Thue equations}.

A couple of other infinite parametric families of Thue equations were
completely solved, see \cite{hh}, \cite{book}, among others 
the parametric family of simplest quartic Thue equations \cite{lp}, \cite{cv} and  
the parametric family of simplest sextic Thue equations \cite{lpv1}, \cite{hoshi2}.

\vspace{0.5cm}

Let $M$ be an algebraic number field with ring of integers $\Z_M$.
Let $F(x,y)\in\Z_M[x,y]$ be  an irreducible binary form of degree 
$n\geq 3$ and let $0\neq \mu\in\Z_M$. As a generalization of Thue equations
consider {\bf relative Thue equations} of type
\[
F(x,y)=\mu\;\;{\rm in}\;\;x,y\in\Z_M.
\]
Using Baker's method S.V.Kotov and V.G.Sprindzuk \cite{ks} gave the first effective upper
bounds for the solutions of relative Thue equations.
Their theorem was extended by several authors.
Applying Baker's method, reduction and enumeration algorithms I.Ga\'al and M.Pohst \cite{gp}
(see also \cite{book}) gave an efficient algorithm for solving relative Thue equations.

Authors considered infinite parametric families of Thue equations
in the relative case, as well. Up to now all these families were
considered over imaginary quadratic fields.
The first of them was the family of
simplest cubic Thue equations, \cite{heuberg1}, \cite{heuberg2}, \cite{heuberg},
\cite{kirsch}. Later on other families of relative Thue equations were also
studied, see e.g. \cite{ziegler1}, \cite{ziegler2}, \cite{jz}.

Let $t$ be an integer parameter,
let $m\geq 1$ be a square-free positive integer, and set $M=\Q(i\sqrt{m})$
with ring of integers $\Z_M$.
In the present paper we consider simplest quartic and simplest sextic Thue equations
in the relative case, over $M$. 
Let 
\[
F_t^{(4)}(x,y)=x^4-tx^3y-6x^2y^2+txy^3+y^4
\]
and let
\[
F_t^{(6)}(x,y)=x^6-2tx^5y-(5t+15)x^4y^2-20x^3y^3+5tx^2
y^4+(2t+6)xy^5+y^6.
\]
We give all solutions of the {\bf infinite parametric families
of simplest quartic and simplest sextic relative Thue equations}. More precisely we
give all solutions of the simplest quartic relative Thue inequalities
\[
|F_t^{(4)}(x,y)|\leq 1 \;\; {\rm in}\;\; x,y\in\Z_M
\]
and of the simplest sextic relative Thue inequalities
\[
|F_t^{(6)}(x,y)|\leq 1 \;\; {\rm in}\;\; x,y\in\Z_M.
\]

\section{Results}

We formulate now our main results. In both Theorems we exclude the
parameters $t\in\Z$ for which the binary form involved is reducible over $\Z$.

\vspace{0.5cm}

\begin{theorem}
Let $t\in\Z$ with $t\neq -3,0,3$.
All solutions of 
\begin{equation}
|F_t^{(4)}(x,y)|\leq 1 \;\; {\rm in}\;\; x,y\in\Z_M
\label{e4}
\end{equation}
are up to sign given by the following:\\
for any $m$ and any $t$: $(x,y)=(0,0),(0,1),(1,0)$,\\
for any $m$ and any $t=1$: $(x,y)= (1,2), (2,-1)$,\\
for any $m$ and any $t=-1$: $(x,y)= (2,1), (-1,2)$,\\
for any $m$ and any $t=4$: $(x,y)=(2,3), (3,-2)$,\\
for any $m$ and any $t=-4$: $(x,y)=(3,2), (-2,3)$,\\
for $m=1$ and any $t$: $(x,y)=(0,i),(i,0)$,\\
for $m=3$ and any $t$: $(x,y)=(\omega,0),(0,\omega),(1-\omega,0),(0,1-\omega)$,\\
for $m=1$ and $t=1$: $(x,y)=(i,2i), (2i,-i)$,\\
for $m=1$ and $t=-1$: $(x,y)= (2i,i), (-i,2i)$,\\
for $m=1$ and $t=4$: $(x,y)= (2i,3i),(3i,-2i)$,\\
for $m=1$ and $t=-4$: $(x,y)= (3i,2i),(-2i,3i)$,\\
for $m=3$ and $t=1$: 
$(x,y)=
(2\omega-2,-\omega+1),
(\omega-1,2\omega-2),
(-2\omega,\omega),
(\omega,2\omega)$,\\
for $m=3$ and $t=-1$: $(x,y)=
(-\omega+1,2\omega-2),
(2\omega-2,\omega-1),
(\omega,-2\omega),
(2\omega,\omega)$,\\
for $m=3$ and $t=4$: $(x,y)=
(3\omega-3,-2\omega+2),
(2\omega-2,3\omega-3),
(2\omega,3\omega),
(3\omega,-2\omega)$,\\
for $m=3$ and $t=-4$: $(x,y)=
(-2\omega+2,3\omega-3),
(3\omega-3,2\omega-2),
(3\omega,2\omega),
(-2\omega,3\omega)$,\\
\label{t4}
where $\omega=(1+i\sqrt{3})/2$.
\end{theorem}

\vspace{0.5cm}

\begin{theorem}
Let $t\in\Z, t\neq -8,-3,0,5$.
All solutions of 
\begin{equation}
|F_t^{(6)}(x,y)|\leq 1 \;\; {\rm in}\;\; x,y\in\Z_M
\label{e6}
\end{equation}
are up to sign given by the following:\\
for any $m$ and any $t$: $(x,y)=(0,0),(0,1),(1,0),(1,-1)$,\\
for $m=1$ and any $t$: $(x,y)=(0, i),( i,0),(i,-i)$,\\
for $m=3$ and any $t$: 
$(x,y)=
(\omega,0),(0,\omega),
(\omega,-\omega),
(1-\omega,0),
(0,\omega-1),
(\omega-1,-\omega+1).
$
\label{t6}
\end{theorem}

\vspace{1cm}

\section{An auxiliary result}

Let $F(x,y)$ be a binary form of degree $n\geq 3$ with rational integer coefficients.
Assume that $f(x)=F(x,1)$ has leading coefficient 1
and distinct real roots $\alpha_1,\ldots,\alpha_n$. 
Let $0<\varepsilon<1,\;\; 0<\eta<1$, and $K\geq 1$. Set
\[
A=\min_{i\not=j}|\alpha_i-\alpha_j|,\;\; B=\min_i\prod_{j\not=i} |\alpha_j-\alpha_i|,
\]
\[
C=\max\left(\frac{K}{(1-\varepsilon)^{n-1}B},1\right),\;
%\]
%\[
C_1=\max\left(\frac{K^{1/n}}{\varepsilon A},\;
(2C)^{1/(n-2)}\right),\;\;
C_2=\max\left(\frac{K^{1/n}}{\varepsilon A},\; 
C^{1/(n-2)}\right),
\]
\[
D=\left(\frac{K}{\eta(1-\varepsilon)^{n-1}AB}\right)^{1/n},\;
E=\frac{(1+\eta)^{n-1}K}{(1-\varepsilon)^{n-1}}.
\]

\noindent
Let $m\geq 1$ be a square-free positive integer, and set $M=\Q(i\sqrt{m})$. 
\\
If $m\equiv 3\; (\bmod \; 4)$, then $x,y\in\Z_M$ can be written as
\[
x=x_1+x_2\frac{1+i\sqrt{m}}{2}=\frac{(2x_1+x_2)+x_2i\sqrt{m}}{2},\;
%\]
%\[
y=y_1+y_2\frac{1+i\sqrt{m}}{2}=\frac{(2y_1+y_2)+y_2i\sqrt{m}}{2}
\]
with $x_1,x_2,y_1,y_2\in\Z$.\\
\noindent
If $m\equiv 1,2\; (\bmod \; 4)$, then $x,y\in\Z_M$ can be written as
\[
x=x_1+x_2i\sqrt{m},\;\; y=y_1+y_2i\sqrt{m}
\]
with $x_1,x_2,y_1,y_2\in\Z$.

Consider the relative Thue inequality
\begin{equation}
|F(x,y)|\leq K \;\; {\rm in}\;\; x,y\in\Z_M.
\label{1}
\end{equation}
We shall use the result of \cite{relthue}:

\begin{lemma}
\label{th1}
Let $(x,y)\in\Z_M^2$ be solutions of (\ref{1}). 
Assume that 
\begin{eqnarray*}
|y|> C_1 & \;\; {\rm if} \;\; &m\equiv 3\; (\bmod \; 4), \\
|y|> C_2 & \;\; {\rm if} \;\; &m\equiv 1,2\; (\bmod \; 4).
\label{f}
\end{eqnarray*}
Then 
\[
x_2y_1=x_1y_2.
\]
\noindent
I. Let $m\equiv 3\; (\bmod \; 4)$.
\begin{eqnarray*}
{\rm IA1.}&\;{\rm If}\;  2y_1+y_2=0,\; {\rm then}\; 2x_1+x_2=0\;{\rm and} \;
|F(x_2,y_2)|\leq \displaystyle{\frac{2^nK}{(\sqrt{m})^n}}.
\label{IA1}\\
{\rm IA2.}&\;{\rm If} \;|2y_1+y_2|\geq 2D,
\; {\rm then} \;
|F(2x_1+x_2,2y_1+y_2)|\leq 2^n E.
\label{IA2}\\
{\rm IB1.}&\;{\rm If}\; y_2=0, \;{\rm then}\; x_2=0\;{\rm and} \;
|F(x_1,y_1)|\leq K.
\label{IB1}\\
{\rm IB2.}&\;{\rm If}\;
|y_2|\geq \displaystyle{\frac{2}{\sqrt{m}} D},\; {\rm then}\;
|F(x_2,y_2)|\leq \frac{2^n}{(\sqrt{m})^n} E.
\label{IB2}
\end{eqnarray*}

\noindent
II. Let $m\equiv 1,2\; (\bmod \; 4)$.
\begin{eqnarray*}
{\rm IIA1.}&\; {\rm If} \;y_1=0, \;{\rm then}\; x_1=0 \;{\rm and} \;
|F(x_2,y_2)|\leq \displaystyle{\frac{K}{(\sqrt{m})^n}}.
\label{IIA1}\\
{\rm IIA2.}&\; {\rm If} \;|y_1|\geq D,\; {\rm then}\;
|F(x_1,y_1)|\leq E.
\label{IIA2}\\
{\rm IIB1.}&\; {\rm If} \; y_2=0,\;{\rm  then}\; x_2=0\;{\rm  and}\;
|F(x_1,y_1)|\leq K.
\label{IIB1}\\
{\rm IIB2.}&\; {\rm If}\; |y_2|\geq \displaystyle{\frac{D}{\sqrt{m}}},\; {\rm then}\;
|F(x_2,y_2)|\leq \frac{E}{(\sqrt{m})^n}.
\label{IIB2}
\end{eqnarray*}

\end{lemma}

\vspace{1cm}

\section{Simplest quartic Thue equations over imaginary quadratic fields}

In this section we turn to the proof of Theorem \ref{t4}.
In our proof we shall use Lemma \ref{th1} and the corresponding results in
the absolute case.

For right hand sides $\pm 1$ J.Chen and P.Voutier \cite{cv} gave all solutions
of simplest quartic Thue equations.

\begin{lemma}
\label{cvlemma}
Let $t\in\Z$ with $t\geq 1,t\neq 3$. All solutions of 
\[
F_t^{(4)}(x,y)=\pm 1 \;\;{\rm in}\;\;x,y\in\Z
\]
are given by 
$(x,y)=(\pm 1,0),(0,\pm 1).$\\ 
Further, for $t=1$ we have
$(x,y)= (1,2),(-1,-2), (2,-1),(-2,1)$\\
and for $t=4$ we have
$(x,y)=(2,3), (-2,-3), (3,-2),(-3,2)$.
\end{lemma} 

For larger right hand sides we can use the statement of 
G.Lettl, A.Peth\H{o} and P.Voutier \cite{lpv}.

\begin{lemma}
\label{l4}
Let $t\in\Z, t\geq 58$ and consider the primitive solutions (i.e. solutions with
$(x,y)=1$) of
\begin{equation}
|F_{t}^{(4)}(x,y)|\leq6t+7\;\;\mathrm{in}\;x,y\in\Z.
\label{f4ri}
\end{equation}
If $(x,y)$ is a solution of (\ref{f4ri}), then every pair in the orbit
\[
\left\{  (x,y),(y,-x),(-x,-y),(-y,x)\right\}
\]
is also a solution. Every orbit has a solution with $y>0,$ $-y\leq x\leq y$.
If an orbit contains a primitive solution, then all solutions in this orbit
are primitive. All solutions of the above inequality with $y>0,-y\leq x\leq y$
are $(0,1),(\pm1,1),(\pm1,2).$
\end{lemma}

\noindent
{\bf Remark 1.} 
Since
\[
F_t^{(4)}(x,y)=F_{-t}^{(4)}(y,x),
\]
it is enough to solve the inequality (\ref{e4}) only for $t> 0$. Also, we have
\[
F_{t}^{(4)}(x,y)=F_{t}^{(4)}(-x,-y)=F_{t}^{(4)}(y,-x)=F_{t}^{(4)}(-y,x).
\]
Therefore, if $(x,y)\in\Z_{M}^{2}$ is solution, then $(y,-x),(-y,x),(-x,-y)$
are solutions, too.

\subsection{Proof of Theorem \ref{t4} for $\bf m\neq1,3$}

\noindent
We use the notation of Lemma \ref{th1}.
Using the estimates of \cite{lpv} for the roots of the polynomial $F_{t}^{(4)}(x,1)$ we 
obtain $A>0.9833$ and $B>58.1$ for $t\geq 58$. Calculating the roots for $0< t<58$
we obtain $A>0.8284, B>4.6114$ for any $t>0,t\neq 3$.
Set
$$\varepsilon=0.1924, \qquad \eta=0.169.$$
For $t>0,t\neq 3$ and a square-free $m$ with $m\neq1,3$
our Lemma \ref{th1} implies:
\begin{corollary}
\label{c6}
Let $(x,y)\in\Z_M^2$ be solutions of (\ref{e4}). 
Assume that 
$$|y|> 6.2741. $$
Then 
\begin{equation*}
x_2y_1=x_1y_2.
\end{equation*}
\noindent
I. Let $m\equiv 3\; (\bmod \; 4)$.
\begin{eqnarray*}
{\rm IA1.}&\;{\rm If}\;  2y_1+y_2=0,\; {\rm then}\; 2x_1+x_2=0\;{\rm and} \;
|F_{t}^{(4)}(x_2,y_2)|\leq 0.326.
\\
{\rm IA2.}&\;{\rm If} \;|2y_1+y_2|\geq 2.618,
\; {\rm then} \;
|F_{t}^{(4)}(2x_1+x_2,2y_1+y_2)|\leq 48.526.
\\
{\rm IB1.}&\;{\rm If}\; y_2=0, \;{\rm then}\; x_2=0\;{\rm and} \;
|F_{t}^{(4)}(x_1,y_1)|\leq 1.
\\
{\rm IB2.}&\;{\rm If}\;
|y_2|\geq 0.989,\; {\rm then}\;
|F_{t}^{(4)}(x_2,y_2)|\leq 0.990.
\end{eqnarray*}

\noindent
II. Let $m\equiv 1,2\; (\bmod \; 4)$.
\begin{eqnarray*}
{\rm IIA1.}&\; {\rm If} \;y_1=0, \;{\rm then}\; x_1=0 \;{\rm and} \;
|F_{t}^{(4)}(x_2,y_2)|\leq 0.25.
\\
{\rm IIA2.}&\; {\rm If} \;|y_1|\geq 1.309,\; {\rm then}\;
|F_{t}^{(4)}(x_1,y_1)|\leq 3.032.
\\
{\rm IIB1.}&\; {\rm If} \; y_2=0,\;{\rm  then}\; x_2=0\;{\rm  and}\;
|F_{t}^{(4)}(x_1,y_1)|\leq 1.
\\
{\rm IIB2.}&\; {\rm If}\; |y_2|\geq 0.925,\; {\rm then}\;
|F_{t}^{(4)}(x_2,y_2)|\leq 0.7582.
\end{eqnarray*}
\end{corollary}

\noindent
{\bf Case I}. $m\equiv 3\; (\bmod \; 4)$

a) Assume that $|y|> 6.2741$. By Corollary \ref{c6} we have\\
If $y_2=0$, then by IB1 we have $x_2=0$ and 
using Lemma \ref{cvlemma} $|F_{t}^{(4)}(x_1,y_1)|\leq 1$ implies $|x_1|,|y_1|\leq 3$.
This contradicts to $|y|> 6.2741$. 

If $y_2\neq 0$, then IB2 implies $|F_{t}^{(4)}(x_2,y_2)|\leq 0.990$, whence 
$y_2=0$, a contradiction. 

Therefore $|y|> 6.2741$ is not possible.

b) Consider now 
$|y|\leq  6.2741$. By Remark 1 if $(x,y)$ is a solution then 
so also is $(y,-x)$. 
As above we obtain that $|x|>6.2741 $ is not possible, 
hence $|x|\leq 6.2741 $.

We enumerate all $x,y$ with $|x|\leq 6.2741 $ and $|y|\leq  6.2741$ and we
obtain the solutions 
$$(x,y)=(0,0),(0,\pm1),(\pm1,0).$$
Additionally we have up to sign\\
for $t=1$: $(x,y)=(1,2), (2,-1)$  and\\
for $t=4$:  $(x,y)=(2,3),(3,-2)$.\\

\vspace{0.5cm}

\noindent
{\bf Case II}. $m\equiv 1,2\; (\bmod \; 4)$\\
Similar to Case I, we obtain the same solutions. \\

\noindent
According to Remark 1 we proved Theorem \ref{t4} for all $t$ with
$t\neq -3,0-3$ and for $m\neq 1,3$.

\subsection{Proof of Theorem \ref{t4} for $\bf m=1$}

Set
$$\varepsilon=0.1792, \qquad \eta=0.0308.$$
For $t>0, t\neq 3$ and $m=1$ Lemma \ref{th1} implies:
\begin{corollary}
Let $(x,y)\in\Z_M^2$ be a solution of (\ref{e4}). Assume
\begin{eqnarray*}
|y|> 6.736.
\end{eqnarray*}
Then 
\begin{equation*}
x_2y_1=x_1y_2.
\end{equation*}
Further,
\begin{eqnarray*}
{\rm IIA1.}&\; {\rm if} \;y_1=0, \;{\rm then}\; x_1=0 \;{\rm and} \;
|F_{t}^{(4)}(x_2,y_2)|\leq 1,
\\
{\rm IIA2.}&\; {\rm if} \;|y_1|\geq 1.98,\; {\rm then}\;
|F_{t}^{(4)}(x_1,y_1)|\leq 1.981,
\\
{\rm IIB1.}&\; {\rm if} \; y_2=0,\;{\rm  then}\; x_2=0\;{\rm  and}\;
|F_{t}^{(4)}(x_1,y_1)|\leq 1,
\\
{\rm IIB2.}&\; {\rm if}\; |y_2|\geq 1.98,\; {\rm then}\;
|F_{t}^{(4)}(x_2,y_2)|\leq 1.981.
\end{eqnarray*}\\
\end{corollary}

a) Assume $|y|> 6.736$. By the above Corollary we deduce:

If $y_1=0$, then by IIA1 we have $|F_{t}^{(4)}(x_2,y_2)|\leq 1$, whence by Lemma \ref{cvlemma}
$|y_2|\leq 3$, contradicting 
$|y|> 6.736$.

If $|y_1|>3$, then by IIA2 we have $|F_{t}^{(4)}(x_1,y_1)|\leq 1$, whence by Lemma \ref{cvlemma}
$|y_1|\leq 3$, a contradiction.

Therefore only $|y_1|=1,2,3$ is possible.

Using IIB1 and IIB2 we similarly obtain that only $|y_2|=1,2,3$ is possible.
But $|y_1|=1,2,3$, $|y_2|=1,2,3$ contradicts $|y|> 6.736$.

b) Hence only $|y|\leq 6.736$ is possible. If $(x,y)$ is a solution, then so also is
$(y,-x)$ therefore we must also have $|x|\leq 6.736$. Enumerating the set 
$(x,y)\in\Z_M^2$ with
$|x|,|y|\leq 6.736$ we obtain 
$$(x,y)=(0,0),(0,\pm1),(\pm1,0),(0,\pm i),(\pm i,0).$$

\noindent
Additionally we have up to sign\\
for $t=1$ $(x,y)=(1,2),  (i,2i), (2,-1), (2i,-i)$ and\\
for $t=4$ $(x,y)=(2,3),  (2i,3i), (3,-2), (3i,-2i)$.\\

\noindent
According to Remark 1 we have proved Theorem \ref{t4} for all $t$ with
$t\neq -3,0,3$ and for $m=1$.

\subsection{Proof of Theorem \ref{t4} for $\bf m=3$}

{\bf First we assume   $t\geq 58$}. Then  $A>0.9833$, and $B> 58.1$. Set
$$\varepsilon=0.6273, \qquad \eta=0.0361.$$

\begin{corollary}
Let $(x,y)eZ_M^2$ be solutions of (\ref{e4}) and let $m=3$. 
Assume $t\geq 58$ and
$$|y|> 1.621.$$
Then 
\begin{equation*}
x_2y_1=x_1y_2.
\end{equation*}
\noindent
Further
\begin{eqnarray*}
{\rm IA1.}&\;{\rm if}\;  2y_1+y_2=0,\; {\rm then}\; 2x_1+x_2=0\;{\rm and} \;
|F_{t}^{(4)}(x_2,y_2)|\leq 1.778,
\\
{\rm IA2.}&\;{\rm if} \;|2y_1+y_2|\geq 3.497,
\; {\rm then} \;
|F_{t}^{(4)}(2x_1+x_2,2y_1+y_2)|\leq 343.753,
\\
{\rm IB1.}&\;{\rm if}\; y_2=0, \;{\rm then}\; x_2=0\;{\rm and} \;
|F_{t}^{(4)}(x_1,y_1)|\leq 1,
\\
{\rm IB2.}&\;{\rm if}\;
|y_2|\geq 2.019,\; {\rm then}\;
|F_{t}^{(4)}(x_2,y_2)|\leq 38.195.
\end{eqnarray*}
\end{corollary}

a)Assume $|y|> 1.621$. Then by the above Corollary:

If $2y_1+y_2=0$, then by IA1 $2x_1+x_2=0$ and $|F_{t}^{(4)}(x_2,y_2)|\leq 1.778$.
By Lemma \ref{cvlemma}
this later inequality implies $(x_2,y_2)=(0,0),(0,\pm1),(\pm1,0)$. However
for $(x_2,y_2)=(0,\pm1),(\pm1,0)$ one of the
the equations $2y_1+y_2=0$ and $2x_1+x_2=0$
have no integer solutions in $y_1$, resp. $x_1$. If $(x_2,y_2)=(0,0)$ then 
$2y_1+y_2=0$ implies $y_1=0$, but $(y_1,y_2)=(0,0)$ contradicts $|y|> 1.621$.

If $|2y_1+y_2|>3.497$, then IA2 implies $|F_{t}^{(4)}(2x_1+x_2,2y_1+y_2)|\leq 343.753$.
Using Lemma \ref{l4} we can easily list all primitive and non-primitive solutions
of this inequality and we always have $|2y_1+y_2|\leq 4$.

Therefore only $|2y_1+y_2|=1,2,3,4$ is possible.

Using IB1 and IB2 we similarly obtain that only $|y_2|=1,2$ is possible.
The equations $|2y_1+y_2|=1,2,3,4$, $|y_2|=1,2$ leave only a few possible values for
$(y_1,y_2)$.

b) If $|x|> 1.621$, then we similarly obtain $|2x_1+x_2|=1,2,3,4$, $|x_2|=1,2$, since
if $(x,y)$ is a solution, then so also is $(y,-x)$.

c1) If $|x|> 1.621$ and $|y|> 1.621$ then we test the finite set
$|2x_1+x_2|=1,2,3,4$, $|x_2|=1,2$, $|2y_1+y_2|=1,2,3$, $|y_2|=1,2$.

c2) If $|x|> 1.621$ and $|y|\leq 1.621$ then we test the finite set 
$|2x_1+x_2|=1,2,3,4$, $|x_2|=1,2$, $|y|\leq 1.621$.

c3) If $|x|\leq 1.621$ and $|y|> 1.621$ then we test the finite set 
$|x|\leq 1.621$, $|2y_1+y_2|=1,2,3,4$, $|y_2|=1,2$.

c4) Finally, if $|x|\leq 1.621$ and $|y|\leq 1.621$ then we test this finite set.\\

\noindent
All together up to sign we get the following solutions for arbitrary $t\geq 58$:\\
$(x,y)=(0,0),(1,0),(0,1),(\omega,0),(0,\omega),(1-\omega,0),(0,1-\omega)$.\\

\vspace{1cm}

\noindent
{\bf Let now $0<t< 58$.}
Considering the roots of the polynomial $F_{t}^{(4)}(x,1)=0$ for these parameters 
we obtain $A> 0.8284$, $B> 4.6114$. Set
$$\varepsilon=0.0348, \qquad \eta=0.0005.$$

\begin{corollary}
Let $m=3$ and $0<t< 58$. Let $(x,y)\in\Z_M^2$ be a solution of (\ref{e4})
and assume
$$|y|> 34.688.$$
Then 
\begin{equation*}
x_2y_1=x_1y_2.
\end{equation*}
Further
\begin{eqnarray*}
{\rm IA1.}&\;{\rm if}\;  2y_1+y_2=0,\; {\rm then}\; 2x_1+x_2=0\;{\rm and} \;
|F_{t}^{(4)}(x_2,y_2)|\leq 1.778,
\\
{\rm IA2.}&\;{\rm if} \;|2y_1+y_2|\geq 9.824,
\; {\rm then} \;
|F_{t}^{(4)}(2x_1+x_2,2y_1+y_2)|\leq 17.825,
\\
{\rm IB1.}&\;{\rm if}\; y_2=0, \;{\rm then}\; x_2=0\;{\rm and} 
|F_{t}^{(4)}(x_1,y_1)|\leq 1,
\\
{\rm IB2.}&\;{\rm if}\;
|y_2|\geq 5.672,\; {\rm then}\;
|F_{t}^{(4)}(x_2,y_2)|\leq 1.981.
\end{eqnarray*}
\end{corollary}

\noindent
a) Assume $|y|> 34.688$. Then by the above Corollary we have:

If $y_2=0$, then by IB1 $|F_{t}^{(4)}(x_1,y_1)|\leq 1$. 
By Lemma \ref{cvlemma} we know the possible solutions $y_1$.
These, together with $y_2=0$ contradict $|y|> 34.688$.

If  $|y_2|\geq 5.672$, then by IB2 $|F_{t}^{(4)}(x_2,y_2)|\leq 1$, whence 
by Lemma \ref{cvlemma} $|y_2|\leq 3$, a contradiction. 

Therefore only $|y_2|=1,2,3,4,5$ is possible.

If $2y_1+y_2=0$, then by IA1 we have $2x_1+x_2=0$ and $|F_{t}^{(4)}(x_2,y_2)|\leq 1$.
From Lemma \ref{cvlemma} we get the possible values of $y_2$ and
we calculate $y_1$ from $2y_1+y_2=0$. These are in contradiction with $|y|> 34.688$.

If $|2y_1+y_2|\geq 9.824$ then by IA2 we have $|F_{t}^{(4)}(2x_1+x_2,2y_1+y_2)|\leq 17$.
Using Magma \cite{magma} we solve the equation $F_{t}^{(4)}(2x_1+x_2,2y_1+y_2)=d$ for 
all $t\leq 58$ and $|d|\leq 17$ and list the solutions. All these solutions
contradict
$|2y_1+y_2|>9.824$.

Therefore only $|2y_1+y_2|=1,\ldots,9$ is possible.

In the set $|y_2|=1,2,3,4,5$, $|2y_1+y_2|=1,\ldots,9$ all $y=y_1+\omega y_2$ have
absolute values less than $34.688$ which is in contradiction with $|y|> 34.688$.

\noindent
b) Therefore only $|y|\leq  34.688$ is possible. 
If $(x,y)$ is a solution, then so also is $(y,-x)$ therefore we
similarly obtain $|x|\leq  34.688$.
Enumerating all $x,y$ with these properties we obtain up to sign the following solutions:\\
for arbitrary $t$: $(1,0),(0,1),(\omega,0),(0,\omega),(1-\omega,0),(0,1-\omega)$,\\
for $t=1$: $
(1,2),
(2,-1),
(2\omega-2,-\omega+1),
(\omega-1,2\omega-2),
(-2\omega,\omega),
(\omega,2\omega),$\\
for $t=4$: $
(2,3),(3,-2),
(3\omega-3,-2\omega+2),
(2\omega-2,3\omega-3),
(2\omega,3\omega),
(3\omega,-2\omega).$\\

\noindent
According to Remark 1 we have proved Theorem \ref{t4}
for all $t$ with $t\neq -3,0,3$ and $m=3$.

\vspace{1cm}

\section{Simplest sextic Thue equations over imaginary quadratic fields}

In this section we turn to the proof of Theorem \ref{t6}. In our proof we shall use 
Lemma \ref{th1} and the corresponding results in the absolute case.

G.Lettl, A.Peth\H o, and P.Voutier \cite{lpv1}
and
A.Hoshi \cite{hoshi2} gave all solutions in rational integers 
of the equation $F_t^{(6)}(x,y)=\pm 1 $ 
for all parameters.

\begin{lemma}
\label{l61}
Let $t\in\Z,t\neq -8,-3,0,5$.
All solutions of 
\[
F_t^{(6)}(x,y)=\pm 1 \;\;{\rm in}\;\;x,y\in\Z
\]
are given by 
\[
(x,y)=(\pm 1,0),(0,\pm 1),(1,-1),(-1,1).
\]
\end{lemma}

For larger right hand sides we shall use the statement of 
G.Lettl, A.Peth\H{o} and P.Voutier \cite{lpv}.

\begin{lemma}
\label{l6}
Let $t\in\Z,t\geq89$ and consider the primitive solutions of
\begin{equation}
|F_t^{(6)}(x,y)|\leq120t+323\;\;{\rm in}\;x,y\in\Z.
\label{f6ri}
\end{equation}
If $(x,y)$ is a solution of the above inequality, then every pair in the orbit
\[
\left\{  (x,y),(-y,x+y),(-x-y,x),(-x,-y),(y,-x-y),(x+y,-x)\right\}
\]
is also a solution. Every orbit has a solution with $y>0,$ $-y/2<x\leq y$. If
an orbit contains one primitive solution, then all solutions in this orbit are
primitive. All solutions of the above inequality with $y>0,-y/2<x\leq y$ are
$(0,1),(1,1),(1,2),(-1,3)$.
\end{lemma}

\noindent
{\bf Remark 2.} 
Since
\[
F_t^{(6)}(x,y)=F_{-t-3}^{(6)}(y,x),
\]
it is enough to solve the inequality (\ref{e6}) only for $t\geq -1,t\neq 0,5$. Also, we have
\[
F_{t}^{(6)}(x,y)=F_{t}^{(6)}(-y,x+y)=F_{t}^{(6)}(-x-y,x)=F_{t}^{(6)}(-x,-y)
=F_{t}^{(6)}(y,-x-y)=F_{t}^{(6)}(x+y,-x).
\]
Therefore, if $(x,y)\in\Z_{M}^{2}$ is solution, then
$(-y,x+y),(-x-y,x),(-x,-y),(y,-x-y),(x+y,-x)$ are solutions, as well.

\subsection{Proof of Thereom \ref{t6} for $m\neq 1,3$}

Using the estimates of \cite{lpv} for $t\geq 89$, we obtain 
$A>0.4986$ and $B>101.83$. Calculating the roots of $F_t^{(6)}(x,1)$ for $-1\leq t< 89, t\neq 0,5$
finally we get $A>0.4646$, $B>3.3121$. Set
$$\varepsilon=0.12, \qquad \eta=0.23.$$

Our Lemma \ref{th1} implies for $m\neq1,3$:
\begin{corollary}
\label{c12}
Let $t\geq -1,t\neq 0,5$ and let $(x,y)\in\Z_M^2$ be solutions of (\ref{e6}). 
Assume that 
$$|y|> 17.937. $$
Then 
\begin{equation*}
x_2y_1=x_1y_2.
\end{equation*}
\noindent
I. Let $m\equiv 3\; (\bmod \; 4)$.
\begin{eqnarray*}
{\rm IA1.}&\;{\rm If}\;  2y_1+y_2=0,\; {\rm then}\; 2x_1+x_2=0\;{\rm and} \;
|F_t^{(6)}(x_2,y_2)|\leq 0.1866.
\\
{\rm IA2.}&\;{\rm If} \;|2y_1+y_2|\geq 2.6453,
\; {\rm then} \;
|F_t^{(6)}(2x_1+x_2,2y_1+y_2)|\leq 341.42.
\\
{\rm IB1.}&\;{\rm If}\; y_2=0, \;{\rm then}\; x_2=0\;{\rm and} 
|F_t^{(6)}(x_1,y_1)|\leq 1.
\\
{\rm IB2.}&\;{\rm If}\;
|y_2|\geq 0.99983,\; {\rm then}\;
|F_t^{(6)}(x_2,y_2)|\leq 0.9954.
\end{eqnarray*}

\noindent
II. Let $m\equiv 1,2\; (\bmod \; 4)$.
\begin{eqnarray*}
{\rm IIA1.}&\; {\rm If} \;y_1=0, \;{\rm then}\; x_1=0 \;{\rm and} \;
|F_t^{(6)}(x_2,y_2)|\leq 0.125.
\\
{\rm IIA2.}&\; {\rm If} \;|y_1|\geq 1.3227,\; {\rm then}\;
|F_t^{(6)}(x_1,y_1)|\leq 5.3347.
\\
{\rm IIB1.}&\; {\rm If} \; y_2=0,\;{\rm  then}\; x_2=0\;{\rm  and}\;
|F_t^{(6)}(x_1,y_1)|\leq 1.
\\
{\rm IIB2.}&\; {\rm If}\; |y_2|\geq 0.93526,\; {\rm then}\;
|F_t^{(6)}(x_2,y_2)|\leq 0.66684.
\end{eqnarray*}
\end{corollary}

\noindent
{\bf Case I}.  $m\equiv 3\; (\bmod \; 4)$ \\
a) Assume that $|y|> 17.937$. Then by Corollary \ref{c12} we have:\\
If $y_2=0$, then by IB1 we have $|F_t^{(6)}(x_1,y_1)|\leq 1$.
By Lemma \ref{l61}
the possible values of $y_1$ and $y_2=0$ contradict $|y|> 17.937$.\\
If $|y_2|\geq 1$ then by IB2 we have $|F_t^{(6)}(x_2,y_2)|=0$, whence $y_2=0$,
a contradiction.\\
b) Therefore only $|y|\leq 17.937$ is possible. If $(x,y)$ is a solution, then by 
Lemma \ref{l6} $(-x-y,x)$ is also a solution. Hence we also must have 
$|x|\leq 17.937$.

\noindent
{\bf Case II}. $m\equiv 1,2\; (\bmod \; 4)$\\
We similarly obtain that only $|x|\leq 17.9365$, $|y|\leq 17.9365$
is possible.
\\

\noindent
In both cases we enumerate all possible 
solutions with $|x|\leq 17.9365$ and $|y|\leq 17.9365$
and finally up to sign obtain
$$(x,y)=(0,0),(0,1),(1,0),(1,-1)$$
According to Remark 2 these are solutions for $m\neq 1,3$
and all $t$ with $t\neq -8,-3,0,5$.

\subsection{Proof of Thereom \ref{t6} for $\bf m=1$}

We set
$$\varepsilon=0.11, \qquad \eta=0.02.$$
For $t>-1,t\neq 0,5$ and $m=1$ Lemma \ref{th1} implies:

\begin{corollary}
Let $(x,y)\in\Z_M^2$ be a solution (\ref{e6}). Assume
\begin{equation*}
|y|> 19.5671.
\end{equation*}
Then 
\begin{equation*}
x_2y_1=x_1y_2.
\end{equation*}
Further, 
\begin{eqnarray*}
{\rm IIA1.}&\; {\rm if} \;y_1=0, \;{\rm then}\; x_1=0 \;{\rm and} \;
|F_t^{(6)}(x_2,y_2)|\leq 1,
\\
{\rm IIA2.}&\; {\rm if} \;|y_1|\geq 1.9685,\; {\rm then}\;
|F_t^{(6)}(x_1,y_1)|\leq 1.9772,
\\
{\rm IIB1.}&\; {\rm if} \; y_2=0,\;{\rm  then}\; x_2=0\;{\rm  and}\;
|F_t^{(6)}(x_1,y_1)|\leq 1,
\\
{\rm IIB2.}&\; {\rm if}\; |y_2|\geq 1.9865,\; {\rm then}\;
|F_t^{(6)}(x_2,y_2)|\leq 1.9772.
\end{eqnarray*}
\end{corollary}

a) Assume $|y|> 19.5671$.\\

If $y_2=0$, then by IIB1 of the above Corollary 
we have $|F_t^{(6)}(x_1,y_1)|\leq 1$. 
By Lemma \ref{l61} the possible values of 
$y_1$ and $y_2=0$ contradict $|y|> 19.5671$.\\

If $|y_2|\geq 2$, then by IIB2 $|F_t^{(6)}(x_2,y_2)|\leq 1$ which implies
by Lemma \ref{l61} $|y_2|\leq 1$,
a contradiction.

Therefore only $|y_2|=1$ is possible.

IIA1 and IIA2 similarly implies that only $|y_1|=1$ is possible.
But $|y_1|=1$, $|y_2|=1$ contradict $|y|> 19.5671$.

b) Therefore only $|y|\leq 19.5671$ is possible. 
If $(x,y)$ is a solution, then by Lemma \ref{l6} so also is $(-x-y,x)$, hence 
we also have $|x|\leq 19.5671$. 
Enumerating all possible values of $x,y$ 
with $|x|\leq 19.5671$, $|y|\leq 19.5671$
up to sign we get the solutions
$$(x,y)=(0,0),(0,1),(1,0),(1,-1),(0, i),( i,0),(i,-i).$$

\noindent
According to Remark 2 these are all solutions for $m=1$ and all $t$ with
$t\neq -8,-3,0,5$.

\subsection{Proof of Theorem \ref{t6} for $\bf m=3$}

{\bf Assume $t\geq 89$}. Then we have $A>0.4986$ and $B>101.83$. 
Set
$$\varepsilon=0.41, \qquad \eta=0.02\; .$$

\begin{corollary}
Let $(x,y)\in\Z_M^2$ be a solution of (\ref{e6}). Assume that 
$$|y|> 4.8917.$$
Then 
\begin{equation*}
x_2y_1=x_1y_2.
\end{equation*}
\noindent
Further,
\begin{eqnarray*}
{\rm IA1.}&\;{\rm if}\;  2y_1+y_2=0,\; {\rm then}\; 2x_1+x_2=0\;{\rm and} \;
|F_t^{(6)}(x_2,y_2)|\leq 2.3703,
\\
{\rm IA2.}&\;{\rm if} \;|2y_1+y_2|\geq 3.0965,
\; {\rm then} \;
|F_t^{(6)}(2x_1+x_2,2y_1+y_2)|\leq 988.372,
\\
{\rm IB1.}&\;{\rm if}\; y_2=0, \;{\rm then}\; x_2=0\;{\rm and} 
|F_t^{(6)}(x_1,y_1)|\leq 1,
\\
{\rm IB2.}&\;{\rm if}\;
|y_2|\geq 1.7877,\; {\rm then}\;
|F_t^{(6)}(x_2,y_2)|\leq 36.606.
\end{eqnarray*}
\end{corollary}

\noindent
a) Assume $|y|> 4.8917$.\\
If $2y_1+y_2=0$, then by IA1 $2x_1+x_2=0$ and $|F_t^{(6)}(x_2,y_2)|\leq 2$.
By Lemma \ref{l6} we have the solutions of this inequality. 
Only $x=0,y=0$ is possible, contradicting $|y|> 4.8917$. \\
If $|2y_1+y_2|\geq 4$ then by IA2 $|F_t^{(6)}(2x_1+x_2,2y_1+y_2)|\leq 988.372$.
Considering the possible primitive and non-primitive solutions of
this inequality, Lemma \ref{l6} implies $|2y_1+y_2|\leq 3$.\\
Therefore only $|2y_1+y_2|=1,2,3$ is possible.\\
If $y_2=0$, then by IB1 $x_2=0$. The possible values of $x_1,y_1$
we obtain from $|F_t^{(6)}(x_1,y_1)|\leq 1$
by Lemma \ref{l61}. These contradict $|y|> 4.8917$.\\
If $|y_2|\geq 2$, then by IB2 $|F_t^{(6)}(x_2,y_2)|\leq 36$. 
Using Lemma \ref{l6} we consider the primitive and
non-primitive solutions of this inequality and we obtain $|y_2|\leq 2$.\\
Therefore only $|y_2|=1,2$ is possible.

b) If $|x|> 4.8917$ then we similarly obtain $|2x_1+x_2|=1,2,3$, $|x_2|=1,2$,
 since if $(x,y)$ is a solution, then by Remark 2 so also is $(-x-y,x)$.

c1) If $|x|> 4.8917$, $|y|> 4.8917$ then we test the finite set 
$|2x_1+x_2|=1,2,3$, $|x_2|=1,2$, $|2y_1+y_2|=1,2,3$, $|y_2|=1,2$.

c2) If $|x|> 4.8917$, $|y|\leq  4.8917$ then we test the finite set 
$|2x_1+x_2|=1,2,3$, $|x_2|=1,2$, $|y|\leq  4.8917$.

c3) If $|x|\leq 4.8917$, $|y|>  4.8917$ then we test the finite set 
$|x|\leq 4.8917$, $|2y_1+y_2|=1,2,3$, $|y_2|=1,2$.

c4)  If $|x|\leq 4.8917$, $|y|\leq  4.8917$ then we test this finite set.\\

\noindent
Finally, for $t\geq 89$, $m=3$, all solutions of $|F_t^{(6)}(x,y)|\leq 1$ up to sign are\\
$
(x,y)=(0,0),(1,0),(0,1),(1,-1),
(\omega,0),(0,\omega),
(\omega,-\omega),
(1-\omega,0),
(0,\omega-1),
(\omega-1,-\omega+1).
$\\
According to Remark 2 these are valid for all values of $t$, $t\leq -92$,
as well.

\vspace{1cm}

\noindent
{\bf Assume now $-1\leq t< 89$}. Then we have $A>0.4646$ and $B>3.3121$. 
Set
$$\varepsilon=0.1124, \qquad \eta=0.0195\; .$$

\begin{corollary}
Assume that 
$$|y|> 19.149\; .$$
Then 
\begin{equation*}
x_2y_1=x_1y_2.
\end{equation*}
\noindent
Further,
\begin{eqnarray*}
{\rm IA1.}&\;{\rm if}\;  2y_1+y_2=0,\; {\rm then}\; 2x_1+x_2=0\;{\rm and} \;
|F_t^{(6)}(x_2,y_2)|\leq 2.371,
\\
{\rm IA2.}&\;{\rm if} \;|2y_1+y_2|\geq 3.962,
\; {\rm then} \;
|F_t^{(6)}(2x_1+x_2,2y_1+y_2)|\leq 127.946,
\\
{\rm IB1.}&\;{\rm if}\; y_2=0, \;{\rm then}\; x_2=0\;{\rm and} \;
|F_t^{(6)}(x_1,y_1)|\leq 1,
\\
{\rm IB2.}&\;{\rm if}\;
|y_2|\geq 2.287,\; {\rm then}\;
|F_t^{(6)}(x_2,y_2)|\leq 4.739.
\end{eqnarray*}
\end{corollary}

\noindent
a) Assume $|y|> 19.149$.

If $y_2=0$, then by IB1 $|F_6(x_1,y_1)|\leq 1$.
We consider its solutions by Lemma \ref{l61} and find that 
the possible $y_1,y_2$ are in contradiction with $|y|> 19.149$.

If $|y_2|\geq 3$, then by IB2 $|F_6(x_2,y_2)|\leq 4$. 
Using Magma \cite{magma}
we solve $F_t^{(6)}(x_2,y_2)=d$ for $-1\leq t < 89$ and $|d|\leq 4$
and find that only $y_2=0,\pm 1$ are possible, contradicting $|y_2|\geq 3$.

Hence only $|y_2|=1,2$ is possible. 

b) If $|x|> 19.149$ then similarly we obtain $|x_2|=1,2$,
since if $(x,y)$ is a solution then so also is $(-x-y,x)$.

c1) If $|x|> 19.149$ and $|y|> 19.149$, then we have 
$|x_2|=1,2$ and $|y_2|=1,2$. 
For any possible pair $x_2,y_2$ we parametrize
$x_1,y_1$ with a single parameter, say $z$, using $x_2y_1=x_1y_2$
(e.g. if $x_2=1,y_2=2$, then $x_1=z,y_1=2z$). For $-1\leq t< 89$
and for the possible right hand sides 
we substitute $x_1,x_2,y_1,y_2$ into 
our original equation (\ref{e6}) to determine the parameter $z$.
We do not obtain any solutions this way.

c2) If $|x|\leq 19.149$ and $|y|> 19.149$, then $|y_2|=1,2$ and
we can enumerate all possible $x_1,x_2$.
For all $-1\leq t < 89$ we determine $y_1$ from our original inequality (\ref{e6})
using $x_1,x_2$ and $|y_2|=1,2$. We do not find any solutions.

c3) If $|x|> 19.149$ and $|y|\leq 19.149$, then we proceed similarly.

c4) If $|x|\leq 19.149$ and $|y|\leq 19.149$, then we test this finite set.\\

\noindent
Finally, for $-1\leq t< 89$, $m=3$, all solutions of $|F_t^{(6)}(x,y)|\leq 1$ up to sign are\\
$
(x,y)=(0,0),(1,0),(0,1),(1,-1),
(\omega,0),(0,\omega),
(\omega,-\omega),
(1-\omega,0),
(0,\omega-1),
(\omega-1,-\omega+1).
$
These are valid for all values of $t$, $-92<t\leq -2$, as well.
Therefore we have proved Theorem \ref{t6} for $m=3$.

\section{Computational aspects}

All auxiliary calculations were made by using Maple \cite{maple}.
Testing a great number of possible solutions took a few hours.

The resolution of Thue equations was performed by using Magma \cite{magma}.
In the quartic case we solved 
$F_t^{(4)}(2x_1+x_2,2y_1+y_2)=d$ for 
all $t\leq 58$ and $|d|\leq 17$.
This took a few minutes.
In the sextic case we solved 
$F_t^{(6)}(x_2,y_2)=d$ for $-1\leq t< 89$ and $|d|\leq 4$.
This took about 30 minutes.

\end{document}